\PassOptionsToPackage{table,xcdraw}{xcolor}
\documentclass[journal,twoside,web]{ieeecolor}
\usepackage{lcsys}

\let\labelindent\relax
\usepackage[inline]{enumitem} 

\let\theoremstyle\relax
\usepackage{amsmath,amsthm,amssymb,amsfonts}
\theoremstyle{plain}
\newtheorem{theorem}{Theorem}
\newtheorem{corollary}{Corollary}

\newtheorem{lemma}{Lemma}

\newtheorem{rule_}{Rule}
\theoremstyle{definition}
\newtheorem{definition}{Definition}

\usepackage{xfrac} 

\usepackage{algorithmic}
\usepackage{algorithm}
\usepackage{array}
\usepackage{textcomp}
\usepackage{stfloats}
\usepackage{url}
\usepackage{verbatim}
\usepackage{graphicx}
\usepackage{hyperref}
\usepackage{url}
\usepackage{afterpage}
\usepackage[normalem]{ulem}
\usepackage{colortbl}
\usepackage{balance}
\usepackage{xspace}
\usepackage{xcolor}
\usepackage{multirow}
\usepackage{lscape}
\usepackage{soul}

\usepackage{tikz, pgfplots, pgfplotstable}  
\usetikzlibrary{arrows}                     
\usetikzlibrary{shapes}                     
\usetikzlibrary{positioning}                
\usetikzlibrary{shapes.misc}                
\pgfplotsset{compat=1.13}
\tikzset{cross/.style={cross out,           
                        draw=black,
                        minimum size=2*(#1-\pgflinewidth),
                        inner sep=0pt,
                        outer sep=0pt},
        cross/.default={1pt}}
\tikzset{main node/.style={ellipse,         
                            minimum height=20pt,
                            minimum width=60pt,
                            draw,
                            font=\Large}}

\newcommand{\strat}{\mathcal{H}}
\newcommand{\graph}{\mathcal{G}}
\newcommand{\satis}{\mathcal{S}}
\newcommand{\states}{V}
\newcommand{\transitions}{E}
\newcommand{\constraints}{\Lambda}

\newcommand{\wh}{\texttt{WeaklyHard.jl}} 

\newcommand{\event}{\alpha}
\newcommand{\astring}{\omega}

\newcommand{\vertiii}[1]{{\left\vert\kern-0.25ex\left\vert\kern-0.25ex\left\vert #1
    \right\vert\kern-0.25ex\right\vert\kern-0.25ex\right\vert}}

\newcommand{\ewhc}{EWHC} 

\newcommand{\new}[1]{\color{black}#1\color{black}\xspace}
\newcommand{\lcss}[1]{\color{black}#1\color{black}\xspace}
\newcommand{\removed}[1]{}

\newcommand{\overbar}[1]{\mkern 1.5mu\overline{\mkern-1.5mu#1\mkern-1.5mu}\mkern 1.5mu}

\newcommand{\abs}[1]{\left|#1\right|}
\newcommand{\norm}[1]{\lVert#1\rVert}

\def\Aa{{\mathcal{A}}}
\def\Alifted{{\mathcal{L}}}
\def\R{{\mathbb{R}}}
\def\N{{\mathbb{N}}}
\def\x{x}

\newcommand{\cM}{\texttt{M}}
\newcommand{\cH}{\texttt{H}}
\newcommand{\cR}{\texttt{R}}
\newcommand{\cX}{\texttt{T}}
\newcommand{\cW}{\texttt{X}}

\renewcommand*{\arraystretch}{1.05}

\begin{document}

\title{Stability of Linear Systems under Extended Weakly-Hard Constraints}

\author{Nils Vreman, Paolo Pazzaglia, Victor Magron, Jie Wang, Martina Maggio \IEEEmembership{Senior Member, IEEE}
\thanks{Nils Vreman and Martina Maggio are with the Department of Automatic Control, Lund University, Sweden. Paolo Pazzaglia and Martina Maggio are with the Department of Computer Science, Saarland University, Germany. Victor Magron is withe the Laboratory for Analysis and Architecture of Systems, CNRS, France. Jie Wang is with the Academy of Mathematics and Systems Science, CAS, China}}

\maketitle
\thispagestyle{empty}

\begin{abstract}
Control systems can show robustness to many events, like disturbances and model inaccuracies. It is natural to speculate that they are also robust to sporadic deadline misses when implemented as digital tasks on an embedded platform. This paper proposes a comprehensive stability analysis for control systems subject to deadline misses, leveraging a new formulation to describe the patterns experienced by the control task under different handling strategies. Such analysis brings the assessment of control systems robustness to computational problems one step closer to the controller implementation.
\end{abstract}



\section{Introduction}
\label{sec:intro}
\IEEEPARstart{R}{obustness} is an essential concern in the design of control systems; they must be able to reliably handle nonlinear effects, unmodeled dynamics and noise, as well as delays in signal transmissions and dropped packets.
A lesser known problem concerns the assessment of robustness to \emph{computational issues} when controllers are implemented as periodic tasks in cheap embedded platforms.
Such tasks are expected to execute with real-time guarantees, i.e., their execution must be completed before a well-defined \emph{deadline}, when the control output must be sent to the actuator.
However, it is common in practice~\cite{akesson2020empirical} that tasks do not always complete within their deadline, causing what is called a \emph{deadline miss}.
This may be caused by delays in computation and memory accesses, transient overloads, bugs and other issues.

A popular model to describe real-time systems allowing deadline misses is the \emph{weakly-hard} model~\cite{Bernat:2001}. 
Weakly-hard tasks feature constraints defining a maximum number of deadlines that can be missed (alternatively, a minimum number to be satisfied) in a given number of consecutive periods.
This model is also the focus of this work.
To analyse the effects on the controlled plant, it is necessary to specify also \emph{what happens when the miss is experienced}, both in terms of changes to the control signal and of actions taken to deal with the failed computation~\cite{Pazzaglia:2019}.
An instance that experiences a deadline miss can be allowed to continue executing until completion (and possibly used later), while in other applications it is stopped and discarded instead.

There is however a mismatch between the guarantees that can be obtained for real-time tasks and platforms~\cite{Ernst:2015,choi2019job}, and the analysis available for \emph{control} tasks under the weakly-hard model.
Fewer works deal with \emph{stability} analysis of weakly-hard real-time control tasks, often targeting specific use-cases. 
For instance, the analysis in~\cite{Maggio:2020} is limited to constraints specifying a maximum number of \emph{consecutive} deadline misses.
The results in \cite{Linsenmayer:2017,linsenmayer2020linear}, obtained for {networked \lcss{linear} control systems} having packet dropouts bounded using the weakly-hard model, can not be generalised for \emph{late completions} \lcss{or \emph{sets} of weakly-hard constraints}.
The authors of~\cite{liang2019security,liang2020leveraging} studied safety guarantees of weakly-hard controllers, considering a miss as a discarded computation with a known periodic pattern.
\lcss{%
In \cite{huang2020saw, huang2019formal}, an over-approximation-based approach is proposed to check the safety of nonlinear weakly-hard systems, where misses are treated as discarded computations and the actuator holds its previous value.
Convergence rates (providing sufficient stability guarantees) are analysed in~\cite{Gaukler:2019a}.
A Lyapunov-based stability analysis of nonlinear weakly-hard systems is studied in~\cite{hertneck2021efficient}, with deadline misses treated as packet dropouts.
However, the state-of-the-art listed above lack generalisability to more expressive real-time implementations, such as different deadline miss models or handling strategies.}

This paper aims at filling the gap, by providing a stability analysis that can be applied to a class of generic weakly-hard models and deadline miss handling strategies.
First, we formally extend the weakly-hard model to explicitly consider the strategy used to handle the miss events. 
By leveraging an automaton representation of the sequences allowed by (a set of) extended weakly-hard constraints, we use Kronecker lifting and the joint spectral radius to properly express its stability conditions.
Using the concept of constraint dominance, we prove analytic bounds on the stability of a weakly-hard system with respect to \emph{less dominant} constraints.
Finally, we analyse the stability of the resulting closed-loop systems using {\tt SparseJSR}~\cite{sparsejsr}, which exploits the sparsity pattern that naturally arises in the Kronecker lifted representation.
The proposed analysis calls for modularity and separation of concern, and can be a useful tool to decouple the constraint specification and the control verification.

\section{Background and Notation}
\label{sec:background}
We consider a controllable and fully observable \emph{discrete-time} sampled linear time invariant system, expressed as
\begin{equation}\label{eq:plant}
  P:\quad
  \begin{cases}
  x_{t+1} &= A_p\,x_t + B_p\,u_t \\
  y_{t} &= C_p\,x_t + D_p\,u_t,
  \end{cases}
\end{equation}
where $x_t\in\R^n$, $u_t\in\R^r$ and $y_t\in\R^q$ are the plant state, the control signal and the plant output, sampled at time $t\cdot T$, $T$ is the sampling period, and $t\in\N$.
The plant is controlled by a stabilising, LTI, one-step delay discrete-time controller
\begin{equation}\label{eq:controller}
  C:\quad
  \begin{cases}
  z_{t+1} = A_c\,z_t + B_c\,\left(r_t - y_t\right) \\
  u_{t+1} = C_c\,z_t + D_c\,\left(r_t - y_t\right),
  \end{cases}
\end{equation}
where $z_t\in\R^s$ is the controller's internal state and $r_t\in\R^q$ is the setpoint.
Without loss of generality, we consider $r_t = 0$.

\subsection{Real-time tasks that may miss deadlines}
\label{ssec:whalgebra}

The controller in~\eqref{eq:controller} is implemented as a real-time task $\tau$, and designed to be executed periodically with period $T$ in a real-time embedded platform.
Under nominal conditions the task releases an instance (called \emph{job}) in each period, that should be completed before the release of the next instance.
We denote the sequence of activation instants for $\tau$ with $(a_i)_{i \in \N}$, such that, in nominal conditions, $a_{i+1} = a_i+T$, the sequence of completion instants $(f_i)_{i \in \N}$, and the sequence of job deadlines with $(d_i)_{i \in \N}$, such that $d_i = a_i + T$ (also called \emph{implicit} deadline).
This requirement can be either satisfied or not, leading respectively to deadline hits and misses.

\begin{definition}[Deadline hit and miss]
\label{def:hit} 
    The $i$-th job of a periodic task $\tau$ with period $T$ hits its deadline when $f_i \leq d_i$ and misses its deadline when $f_i > d_i$.
\end{definition}

We refer to both deadline hits and misses using the term \emph{outcome} of a job.
\lcss{Intuitively, each job's outcome is dependent on the characteristics of the remaining tasks executing in the real-time system and the chosen scheduling algorithm.
Given a taskset and a (worst-case) schedule, it is possible to bound the worst-case behaviour of the job outcomes~\cite{Bernat:2001, Ernst:2015}.
This bound is generally denoted using the \emph{weakly-hard model}~\cite{Bernat:2001}.}
\lcss{Following such model,}
    a task $\tau$ may satisfy any combination of these weakly-hard constraints, \lcss{defined as follows.}
    \begin{enumerate*}[label=(\roman*)]
        \item \label{item:mk} $\tau \vdash\overbar{\binom{m}{k}}$: in any window of $k$ consecutive jobs, at most $m$ deadlines are missed;
        \item \label{item:hk} $\tau \vdash\binom{h}{\!\:\!\:k\!\:\!\:}$: in any window of $k$ consecutive jobs, at least $h$ deadlines are hit;
        \item \label{item:cons} $\tau \vdash\overbar{\genfrac{<}{>}{0pt}{}{m}{\!\:\!\:k\!\:\!\:}}$: in any window of $k$ consecutive jobs, at most $m$ \emph{consecutive} deadlines are missed; and
        \item $\tau \vdash\genfrac{<}{>}{0pt}{}{h}{\!\:\!\:k\!\:\!\:}$: in any window of $k$ consecutive jobs, at least $h$ \emph{consecutive} deadlines are hit.
    \end{enumerate*}
    \lcss{In all such cases,} $m,h\in \N$, $k \in \N\setminus \left\{ 0 \right\}$, $m\leq k$, and $h\leq k$.
%
A generic weakly-hard constraint is hereafter denoted with the symbol $\lambda$, while a set of $L$ constraints will be referred to as $\constraints = \{ \lambda_1, \lambda_2, \dots, \lambda_L \}$.

We define a \emph{string} $\astring= \lcss{\langle \alpha_1, \alpha_2, \dots, \alpha_N \rangle}$ as a sequence of $N$ consecutive outcomes, where each outcome $\alpha_i$ is a character in the alphabet $\Sigma = \{\cM, \cH\}$.
We use $\astring \vdash \lambda$ to denote that $\astring$ satisfies the constraint $\lambda$. 
Stating that $\tau \vdash \lambda$ means that all the possible sequences of outcomes that $\tau$ can experience satisfy the corresponding constraint $\lambda$.
The set of such sequences naturally results from the definition of $\lambda$, and is formally defined as the \emph{satisfaction set} as follows~\cite{Bernat:2001}.

\begin{definition}[Satisfaction set $\mathcal{S}_N\left(\lambda\right)$\removed{~of a weakly-hard constraint $\lambda$}] \label{def:satisfaction}
    We denote with $\mathcal{S}_N\left(\lambda\right)$ the set of strings of length $N \geq 1$ that satisfy a constraint $\lambda$.
    Formally, $\mathcal{S}_N\left(\lambda\right) = \{ \astring \in \Sigma^N \mid \astring \vdash \lambda \}$.
\end{definition}
Taking the limit to infinity, the set $\mathcal{S}\left(\lambda\right)$ contains all the strings of infinite length that satisfy $\lambda$.
The notion of \emph{domination} between constraints~\cite{Bernat:2001} then follows.
\begin{definition}[Constraint domination] \label{def:domination}
    Constraint $\lambda_i$ \emph{dominates} $\lambda_j$ (formally, $\lambda_i \preceq \lambda_j$) if $ \mathcal{S}\,( \lambda_i ) \subseteq \mathcal{S} \,( \lambda_j )$. 
\end{definition}

\subsection{Control tasks that may miss deadlines}
\label{sec:back_deadline_miss}
When a control task $\tau$ is implemented on an embedded platform with limited computational power, alongside other applications, it is not uncommon for it to experience deadline misses, even in case of simple control designs (PID, LQG, etc)~\cite{akesson2020empirical,pazzaglia2021adaptive}.
Computational overruns may be caused by, e.g., bursts of interrupts, cache misses, variable execution times of ancillary functions, or other complex interactions.
If such events are rare or temporary, choosing a longer period for the controller to avoid them may result in worse performance and stability margins for nominal conditions~\cite{Pazzaglia:2019}.

Characterising the stability and performance of such controllers requires knowing what happens when a control deadline is missed~\cite{Pazzaglia:2019,Maggio:2020,Vreman:2021}.
In particular, we need a \emph{deadline miss handling strategy} to decide the fate of the job that missed the deadline (and possibly the next ones), and an \emph{actuator mode} to deal with the loss of a new control signal, for example by holding the previous value constant or zeroing it~\cite{schenato09}.
A few handling strategies for periodic controllers have been proposed in literature, the most interesting being \emph{Kill} and \emph{Skip-Next}~\cite{Cervin:2005,Pazzaglia:2019,Maggio:2020}.

\begin{definition}[Kill strategy]
    \label{def:kill} 
    Under the Kill strategy, a job that misses its deadline is terminated immediately. 
    Formally, for the~$i$-th job of $\tau$ either $f_i\leq d_i$ or $f_i=\infty$.
\end{definition}

\begin{definition}[Skip-Next strategy]
    \label{def:skip} 
    Under the Skip-Next strategy, a job that misses its deadline is allowed to continue during the following period.
    Formally, if the $i$-th job of $\tau$ misses its deadline $d_i$, a new deadline $d^{+}_i = d_i + T$ is set for the job, and $a_{i+1} = d^{+}_i$.
\end{definition}

\subsection{Stability analysis techniques based on JSR}
\label{sec:existing}

In~\cite{Maggio:2020}, the authors identify a set of subsequences of hit and missed deadlines, which can be arbitrarily combined to obtain all possible sequences in $\mathcal{S}\,(\overbar{\genfrac{<}{>}{0pt}{}{m}{\!\:\!\:k\!\:\!\:}})$.
The stability analysis of the resulting arbitrary switching system is then obtained by leveraging the {\em Joint Spectral Radius} (JSR)~\cite{rota}.

Given $\ell\in\N\setminus \left\{ 0 \right\} $ and a set of
matrices $\Aa = \{A_1, \ldots, A_\ell \} \subseteq \R^{n\times n}$,
under the hypothesis of arbitrary switching \new{over any sequence $s=\lcss{\langle a_1,a_2,\dots \rangle}$ of indices of matrices in $\Aa$}, the JSR of $\Aa$ is defined by:
\begin{equation}
  \label{jsr} \rho\left(\Aa\right)= \lim_{k\rightarrow\infty}
  \max_{s\in\{1,\ldots,\ell\}^k} \norm{A_{a_k}
  \cdots A_{a_2}A_{a_1}}^{\frac{1}{k}} \,.
\end{equation}
The number $\rho\left(\Aa\right)$ characterizes the maximal asymptotic growth rate of matrix products from $\Aa$ (thus  $\rho(\Aa)<1$ means that the system is asymptotically stable), and is independent of the norm $\norm{\cdot}$ used in \eqref{jsr}.
Existing practical tools such as the JSR Matlab toolbox~\cite{vankeerberghen2014jsr} include multiple algorithms to compute both upper and lower bounds on $\rho\left(\Aa\right)$.

When the switching sequences between the dynamics of $\Aa$ are not arbitrary, but constrained by a graph $\graph$, the so called \emph{constrained joint spectral  radius} (CJSR)~\cite{dai2012gelfand} can be applied.
Introducing ${S}_k(\graph)$ as the set of all possible switching sequences $s$ of length $k$ that satisfy the constraints of a graph $\graph$, the CJSR of $\Aa$ is defined by
\begin{equation}
  \label{cjsr} \rho\left(\Aa,\graph\right)= \lim_{k\rightarrow\infty}
  \max_{s\in{S}_k(\graph)} \norm{A_{a_k}
  \cdots A_{a_2}A_{a_1}}^{\frac{1}{k}} \,.
\end{equation}
In general, computing or approximating the CJSR is harder than using the JSR. 
In~\cite{philippe2016stability}, the authors propose a multinorm-based method to approximate with arbitrary accuracy the CJSR.
Other works~\cite{kozyakin2014berger,xu2020approximation} propose the creation of an arbitrary switching system such that its JSR is equal to the CJSR of the original system, based on a Kronecker lifting method.
This will be also our approach, as detailed later.

In~\cite{parrilo}, the authors propose an efficient approach to compute upper bounds of the JSR based on positive polynomials which can be decomposed as \emph{sums of squares} (SOS). 
%
Finding the coefficients of a polynomial being SOS simplifies to solving an SDP~\cite{lasserre2001global}.
To reduce time and space complexity, a \emph{sparse} variant has been proposed in \cite{sparsejsr} exploiting the sparsity of the input matrices, based on the \emph{term sparsity} SOS (TSSOS) framework~\cite{tssos}.
By contrast, the procedure in~\cite{parrilo} will be denoted hereafter as \emph{dense}.
While providing a more conservative result, the sparse upper bound can be obtained significantly faster if the matrices from $\Aa$ are sparse~\cite{sparsejsr}, e.g., the matrices we analyse in Section~\ref{sec:stability}.

\section{Extended Weakly-Hard Task Model}
\label{sec:model}
To provide a comprehensive analysis framework, we need to examine what occurs in each time interval $(\pi_i)_{i \in \N}$, with $\pi_i = [a_0 + i\cdot T, a_0 + (i+1)\cdot T)$. 
In this context, an extension of the weakly-hard model is required to account for the given deadline miss handling strategy, denoted with the symbol $\strat$.
\begin{definition}[Extended weakly-hard model $\tau \vdash \lambda^{\strat}$]
    \label{def:new-mk}
    A task $\tau$ may satisfy any combination of the four \emph{extended weakly-hard constraints} (\ewhc{}) $\lambda^{\strat}$:
    \begin{enumerate}[label=(\roman*)]
        \item $\tau \vdash \overbar{\binom{m}{k}}^{\strat}$: in any window of $k$ consecutive jobs, at most $m$ intervals lack a job completion;
        \item $\tau \vdash \binom{h}{\!\:\!\:k\!\:\!\:}^{\strat}$:  in any window of $k$ consecutive jobs, at least $h$ intervals have a job completion;
        \item $\tau \vdash \overbar{\genfrac{<}{>}{0pt}{}{m}{\!\:\!\:k\!\:\!\:}}^{\strat}$: in any window of $k$ consecutive jobs, at most $m$ \emph{consecutive} intervals lack a job completion;
        \item $\tau \vdash \genfrac{<}{>}{0pt}{}{h}{\!\:\!\:k\!\:\!\:}^{\strat}$: in any window of $k$ consecutive jobs, at least $h$ \emph{consecutive} intervals have a job completion
    \end{enumerate}
    with $m,h\in \N$, $k \in \N\setminus \left\{ 0 \right\}$, $m\leq k$, and $h\leq k$, while using strategy $\strat$ to handle potential deadline misses.
\end{definition}
The definition above differs from \lcss{the original weakly-hard model of~\cite{Bernat:2001},}
since (i) it explicitly introduces the handling strategy $\strat$; and (ii) it focuses on the presence of a new control command at the end of each time interval $\pi_i$, instead of checking the deadline miss events, which guarantees its applicability also for strategies different than Kill.

We now require an expressive alphabet $\Sigma\left(\strat\right)$ to characterize the behaviour of task $\tau$ in each possible time interval.
For both Kill and Skip-Next strategies, each interval $\pi_i$ contains at most one activated and one completed job.
This restricts the possible behaviours to three cases:
\begin{enumerate}[label=(\roman*)]
    \item a time interval in which the same job is both released and completed is denoted by $\cH$ (\emph{hit});
    \item a time interval in which no job is completed is denoted by $\cM$ (\emph{miss});
    \item a time interval in which no job is released, but a job (released in a previous interval) is completed, is denoted by $\cR$ (\emph{recovery}).
\end{enumerate}
By checking all unique combinations of job activations and completions in each interval, we obtain the alphabets for Kill and Skip-Next as $\Sigma\left(\text{Kill}\right) = \{ \cM, \cH \}$ and $\Sigma\left(\text{Skip-Next}\right) = \{ \cM, \cH, \cR \}$, respectively.
The recovery character $\cR$ is used in the Skip-Next alphabet to identify the late \emph{completion} of a job.
As a consequence, $\cR$ is treated equivalently to $\cH$ when verifying the extended weakly hard constraints (\ewhc{}).

The algebra presented in Section~\ref{ssec:whalgebra} is extended to the new alphabet.
We assign a character of the alphabet $\Sigma\left(\strat\right)$ to each interval $\pi_i$.
A string $\astring = \lcss{\langle \alpha_1,\alpha_2,\dots,\alpha_N \rangle}$ is used to represent a sequence of $N$ outcomes for task $\tau$, with $\alpha_i \in \Sigma\left(\strat\right)$ representing the outcome associated to the interval $\pi_i$. 
%
To enforce only feasible sequences, we introduce an order constraint for the $\cR$ character with the following Rule.
\begin{rule_}
    \label{rule:R}
    For any string $\astring \in \Sigma\left(\text{Skip-Next}\right)^N$, $\cR$ may only directly follow \removed{a miss}$\cM$, or be the initial element of the string.
\end{rule_}

The extended weakly-hard model also inherits all the properties of the original weakly-hard model.
In particular, the satisfaction set of $\lambda^\strat$ can be defined for $N\geq 1$ as $\mathcal{S}_N\,(\lambda^{\strat}) = \{ \astring \in \Sigma\left(\strat\right)^N \mid \astring \vdash \lambda^{\strat} \}$, and the constraint domination still holds as $\lambda^{\strat}_{i} \preceq \lambda^{\strat}_{j}$ if $\mathcal{S}\,( \lambda^{\strat}_{i} ) \subseteq \mathcal{S}\,( \lambda^{\strat}_{j})$.

\subsection{Automaton representation of \ewhc{}}%
Any \ewhc{}, as presented in Definition~\ref{def:new-mk}, can be systematically represented using an \emph{automaton}.
In this paper we build upon the \wh{} automaton model presented in~{\cite{newNilsPaper}}.
Here, a (minimal) automaton $\graph_{\lambda^\strat} = \left(\states_{\lambda^\strat}, \transitions_{\lambda^\strat}\right)$ associated to $\lambda^\strat$ consists of a set of vertices ($\states_{\lambda^\strat}$) and a set of directed labeled edges ($\transitions_{\lambda^\strat}$). 
Each vertex $v_i \in \states_{\lambda^\strat}$ corresponds to a string of outcomes of the extended weakly-hard task executions. 
Trivially, there exists no vertices for strings that do not satisfy the \ewhc{}.
A directed labeled edge $e_{i, j} = \left(v_i, v_j, \event \right) \in \transitions_{\lambda^\strat}$ (also denoted \emph{transition}) connects two vertices iff the outcome $\event \in \Sigma\left(\strat\right)$ -- the edge's label -- appended to the tail vertex's string representation ($v_i$) would result in the string equivalent to the one of the head vertex ($v_j$).
Thus, a random walk in the automaton corresponds to a random string satisfying the \ewhc{}.
In particular, all the walks in the automaton corresponds to \emph{all} strings in $\satis\,(\lambda^{\strat})$.

Since the \wh{} automaton model only uses the binary alphabet $\Sigma = \left\{\cM, \cH\right\}$, we require the additional character $\cR$ to handle the Skip-Next strategy properly.
Recall that both a hit ($\cH$) and a recovery ($\cR$) are considered job completions.
Thus, for the Skip-Next strategy, we post-process the automaton by enforcing that Rule~\ref{rule:R} is honoured and that the corresponding transitions are correct, i.e., switching the labels on some edges from $\cH$ to $\cR$. 
We emphasise that despite the extended automaton model appear similar for the Kill and Skip-Next strategies, the differing transitions of the two automata significantly affect the corresponding closed-loop systems, as will be clear in Section~\ref{sec:stability}.

The \wh{} automaton model also allows for the case where the task $\tau$ is subject to a set of multiple constraints.
Since the stability analysis presented in this paper is invariant to the type (and amount) of the constraints acting on the control task $\tau$, we henceforth say that $\tau$ is subject to a set of \ewhc{} $\constraints^\strat$ (unless stated otherwise).

%
Extracting all transitions in $\transitions_{\constraints^\strat}$ corresponding to a character $\event \in \Sigma\left(\strat\right)$ yields what is generally known as a \emph{directed adjacency matrix}~\cite{xu2012matrix}, denoted here as a \emph{transition matrix}.
\begin{definition}[Transition matrix]
    \label{def:transition}
    Given an automaton $\graph_{\constraints^\strat}$, the \emph{transition matrix} $F_{\event} ( \graph_{\constraints^\strat} ) \in \R^{n_V\times n_V}$, with $n_V =\abs{\states_{\constraints^\strat}} $ and $\event\in\Sigma\left(\strat\right)$, is computed as $F_{\event} ( \graph_{\constraints^\strat} ) = \{f_{i,j}(\event)\}$ with
    \begin{equation*}
        f_{i,j}\left(\event\right)=
        \begin{cases}
            1, &\text{ if } \exists \, e=(v_j,v_i,\event) \in \transitions_{\constraints^\strat} \\
            0, &\text{ otherwise.}
        \end{cases}
        \end{equation*}%
\end{definition}
Since \emph{at most one} successor exists from each vertex with a transition labeled with $\event \in \Sigma\left(\strat\right)$, matrix $F_\event$ will have a column sum of either 1 or 0.
We now introduce a vector $q_t\in \R^{n_V}$ called \emph{G-state}, with $n_V =\abs{\states_{\constraints^\strat}} $, representing the state of the given automaton $\graph_{\constraints^\strat}$ at interval $\pi_t$.
\begin{definition}[G-state $q_t$]\label{def:qt}
    Given an automaton $\graph_{\constraints^\strat}$ and a string $\astring \in \Sigma\left( \strat \right)^N$, $\astring=\lcss{\langle \alpha_1,\alpha_2,\dots,\alpha_N \rangle}$, for $k = \abs{v},\,\, v\in\states_{\constraints^\strat}$, we define $q_t\in \R^{n_V}$, where the $i$-th element $q_{t,i}$ is: 
    \begin{equation*}
        q_{t,i}=
        \begin{cases}
            1, &\text{ if } \lcss{\langle \alpha_{t-k},\dots,\alpha_{t-1} \rangle} \equiv v_i \in \states_{\constraints^\strat} \\
            0, &\text{otherwise}.
        \end{cases}
    \end{equation*}
\end{definition}
The G-state $q_t$ is the vector representation of the vertex \emph{left} at step $t$: here, $q_t=0$ means that the transition at step $t-1$ was infeasible for the automaton. 
Given an arbitrary string $\astring=\lcss{\langle \alpha_1,\dots,\alpha_t,\dots \rangle}$, the G-state dynamics is defined as $q_{t+1} = F_{\event} ( \graph_{\constraints^\strat} )\cdot q_t$, and the following property holds~\cite{xu2012matrix}.
\begin{lemma}[Infeasible sequence]
    \label{cor:Fseqnotinlambda}
    If $\astring \notin \satis_N \,( \constraints^\strat )$, then $F_{\astring} ( \graph_{\constraints^\strat} ) =
    F_{\event_N} ( \graph_{\constraints^\strat} )\cdots F_{\event_2} ( \graph_{\constraints^\strat} )\cdot F_{\event_1} ( \graph_{\constraints^\strat} ) = 0$
\end{lemma}
Thus, if $q_t=0$ for an arbitrary $t$, then $q_{t'}=0$ for $t' \geq t$.

\section{Stability Analysis}
\label{sec:stability}

Using the alphabet $\Sigma\left(\strat\right)$ and the chosen actuator mode (i.e., zeroing, or holding the previous value), we compute the closed-loop behaviour of the controlled system.
We identify one matrix for each dynamics corresponding to an interval $\pi_t$ associated by $\event \in \Sigma\left( \strat \right)$, building the set $\Aa^\strat$.

\textbf{Kill: }%
%
Defining $\tilde x_t^{\,\text{K}} = \left[ x_t^T\,\, z_t^T\,\, u_t^T \right]^T$ as the closed-loop state vector,
we compute the discrete time closed-loop system dynamics $A^{\,\text{K}}_{\cH}$, corresponding to the character $\cH$:
\begin{equation*}
    \tilde x_{t+1}^{\,\text{K}} = A^{\,\text{K}}_{\cH}\,\tilde x_t^{\,\text{K}}, \quad\;
     A^{\,\text{K}}_{\cH} = \begin{bmatrix}
        A_p       & 0    & B_p       \\
        -B_cC_p   & A_c  & -B_cD_p   \\
        -D_cC_p   & C_c  & -D_cD_p   \\
    \end{bmatrix}.
\end{equation*}
For the case of $\cM$, the controller execution terminates prematurely and its states are not updated ($z_{t+1} = z_t$).
Therefore, depending on the actuation mode, the controller output is either zeroed ($u_{t+1} = 0$) or held ($u_{t+1} = u_t$).
The resulting closed-loop system in state-space form is denoted with $A^{\,\text{K}}_{\cM}$:
\begin{equation*}
    \tilde x_{t+1}^{\,\text{K}} = A^{\,\text{K}}_{\cM}\,\tilde x_t^{\,\text{K}}, \quad
   A^{\,\text{K}}_{\cM} = \begin{bmatrix}
        A_p & 0  & B_p \\
        0   & I  & 0   \\
        0   & 0  & \Delta
    \end{bmatrix}.
\end{equation*}
Here, $\Delta = I$ (identity matrix) if the control signal is held and $\Delta = 0$ if zeroed.
The set of dynamic matrices under the Kill strategy is then $\Aa^{\text{Kill}}=\left\{A^{\,\text{K}}_{\cH},A^{\,\text{K}}_{\cM}\right\}$.

\textbf{Skip-Next: }%
For the Skip-Next strategy, we introduce two additional states $\hat x_t$ and $\hat u_t$ storing the old values of $x_t$ and $u_t$ while the controller awaits an update.
The resulting state vector then becomes $\tilde x_t^{\,\text{S}} = \left[ x_t^T\,\,z_t^T\,\, u_t^T\,\, \hat{x}_t^T\,\, \hat{u}_t^T \right]^T$.
When $\pi_t$ is associated to $\cH$, the two additional states mirror the behaviour of the states of which they are storing data.
The resulting closed-loop system is described using $A^{\,\text{S}}_{\cH}$:
\begin{equation*}
    \setlength\arraycolsep{3pt}
    \tilde x_{t+1}^{\,\text{S}} = A^{\,\text{S}}_{\cH}\,\tilde x_t^{\,\text{S}}, \quad
    A^{\,\text{S}}_{\cH} = \begin{bmatrix}
        A_p       & 0    & B_p      & 0 & 0 \\
        -B_cC_p   & A_c  & -B_cD_p  & 0 & 0 \\
        -D_cC_p   & C_c  & -D_cD_p  & 0 & 0 \\
        A_p       & 0    & B_p      & 0 & 0 \\
        -D_cC_p   & C_c  & -D_cD_p  & 0 & 0 \\
    \end{bmatrix}.
\end{equation*}
For the case of $\cM$ in $\pi_t$, $\hat x_t$ and $\hat u_t$ maintain their previous values. The
resulting closed-loop is described by $A^{\,\text{S}}_{\cM}$:
\begin{equation*}
    \setlength\arraycolsep{3pt}
    \tilde x_{t+1}^{\,\text{S}} = A^{\,\text{S}}_{\cM}\,\tilde x_t^{\,\text{S}}, \quad
    A^{\,\text{S}}_{\cM}=\begin{bmatrix}
        A_p & 0  & B_p & 0 & 0 \\
        0   & I  & 0   & 0 & 0 \\
        0   & 0  & \Delta   & 0 & 0 \\
        0   & 0  & 0   & I & 0 \\
        0   & 0  & 0   & 0 & I \\
	\end{bmatrix}.
\end{equation*}
Finally, for the case of $\cR$, the new control command is calculated using the values stored in $\hat x_t$ and $\hat u_t$.
The resulting closed-loop system is described by $A^{\,\text{S}}_{\cR}$:
\begin{equation*}
    \setlength\arraycolsep{3pt}
    \tilde x_{t+1}^{\,\text{S}} = A^{\,\text{S}}_{\cR} \, \tilde x_{t}^{\,\text{S}},\quad
    A^{\,\text{S}}_{\cR} = \begin{bmatrix}
        A_p & 0    & B_p & 0       & 0 \\
        0   & A_c  & 0   & -B_cC_p & -B_cD_p \\
        0   & C_c  & 0   & -D_cC_p & -D_cD_p \\
        A_p & 0    & B_p & 0       & 0 \\
        0   & C_c  & 0   & -D_cC_p & -D_cD_p \\
    \end{bmatrix}.
\end{equation*}
The resulting set of matrices under Skip-Next strategy is then defined as $\Aa^{\text{Skip-Next}}=\left\{A^{\,\text{S}}_{\cH},A^{\,\text{S}}_{\cM},A^{\,\text{S}}_{\cR}\right\}$.

\subsection{Kronecker lifted switching system}%
\label{sec:system_dynamics}

%
Combining the set of system dynamics $\Aa^\strat$ with the associated automaton $\graph_{\constraints^\strat}$, we seek to obtain an equivalent system model based on Kronecker lifting, characterized by a set of matrices denoted by $\Alifted_{\constraints^\strat}$ and behaving as an \emph{arbitrary switching system}, such that $\rho\,(\Alifted_{\constraints^\strat})= \rho\,(\Aa^{\strat},\graph_{\constraints^\strat})$.
In this way, powerful algorithms applicable to arbitrary switching system~\cite{vankeerberghen2014jsr,sparsejsr} can be used to find tight stability bounds.
We build upon the Kronecker lifting approach of~\cite{xu2020approximation}.
Leveraging the vector $q_t$ of Definition~\ref{def:qt}, we introduce the \emph{lifted discrete-time state} $\xi_t\in\mathbb{R}^{n\cdot n_V}$, defined as 
$\xi_t = q_t\otimes x_t$, where $n_V = \abs{\states_{\constraints^\strat}}$ and $\otimes$ is the Kronecker product.
By construction, $\xi_t$ is a vector composed of $n_V$ blocks of size $n$, where at most one block is equal to $x_t$ and all other blocks are equal to the $0$ vector.
Then, we build a set of lifted matrices $P_{\event} ( \graph_{\constraints^\strat} )\in\mathbb{R}^{n\cdot n_V\times n\cdot n_V}$, which incorporates both the system dynamics and the possible transitions given a certain outcome $\event\in\Sigma\left(\strat\right)$:
\begin{equation}\label{eq:lifted_matrix}
    P_{\event} ( \graph_{\constraints^\strat} ) = F_{\event} ( \graph_{\constraints^\strat} )\otimes A_\event^\strat,\quad \event \in \Sigma \left( \strat \right).
\end{equation}
The lifted dynamics of the closed loop system then become 
$
    \xi_{t+1} = P_{\event} ( \graph_{\constraints^\strat} )\cdot\xi_t.
$
Formally, we obtain a system composed of a set of switching dynamic matrices, $\Alifted_{\constraints^\strat}$.
\begin{definition}[Lifted switching set $\Alifted_{\constraints^\strat}$]
    \label{def:switching_set}
    Given a set of dynamic matrices $\Aa^{\strat}$ and an automaton $\graph_{\constraints^\strat}$, the switching set $\Alifted_{\constraints^\strat}$ is defined as:
    $
    \Alifted_{\constraints^\strat} = \left\{ P_{\event} ( \graph_{\constraints^\strat} ) \,\, | \,\, \event \in \Sigma\left(\strat\right) \right\}. $
\end{definition}%
Leveraging the mixed-product property of $\otimes$ and introducing a proper submultiplicative norm, it is possible to prove that $\rho\left(\Alifted_{\constraints^\strat}\right)= \rho\,(\Aa^{\strat},\graph_{\constraints^\strat})$.
For more details and a formal proof we refer the interested reader to~\cite{xu2020approximation}.

\subsection{Extended weakly hard and JSR properties}
\label{sec:analytic_results}
We now provide a general relation between \emph{all} \ewhc{}s in terms of the joint spectral radii.
%
\begin{theorem}[JSR dominance]
    \label{th:rho_dominance_general}
    Given $\lambda_1^\strat$ and $\lambda_2^\strat$ as arbitrary \ewhc{}s, if $\lambda_2^\strat \preceq \lambda_1^\strat$ then
    $
        \rho\bigl(\Alifted_{\lambda_2^\strat}\bigr) \leq \rho\bigl(\Alifted_{\lambda_1^\strat}\bigr).
    $
    \begin{proof}
        From Equation~\eqref{jsr}, for a generic \ewhc{} $\lambda^\strat$,
        \begin{equation*}
            \rho\left(\Alifted_{\lambda^\strat}\right) = \lim_{\ell\rightarrow
            \infty}\rho_\ell\left(\Alifted_{\lambda^\strat}\right), \;\, \rho_\ell\left(\Alifted_{\lambda^\strat}\right) =
            \max_{a \in
            \satis_\ell\left(\lambda^\strat\right)}\norm{A_{a}}^{1/\ell}.
        \end{equation*}
        Definition~\ref{def:domination} gave us that $\lambda^\strat_2 \preceq \lambda^\strat_1$ iff $\satis(\lambda^\strat_2) \subseteq \satis(\lambda^\strat_1)$.
        Thus, if for a string $b$ it holds that $b \in \satis_\ell \left( \lambda^\strat_2 \right)$, then it also holds that $b \in \satis_\ell \left( \lambda^\strat_1 \right)$.
        The set of all possible $A_{b}$ is thus included in the set of all possible $A_{a},\, a \in \satis_\ell \left( \lambda^\strat_1 \right)$, thus:
        \begin{equation*}
            \max_{b \in
            \satis_\ell(\lambda^\strat_2)}\norm{A_{b}}^{1/\ell} \leq
            \max_{a \in
            \satis_\ell(\lambda^\strat_1)}\norm{A_{a}}^{1/\ell}, \quad
            \forall \ell\in\mathbb{N}\setminus 0.
        \end{equation*}
        The theorem follows immediately when $\ell\rightarrow \infty$.
    \end{proof}
\end{theorem}
%
Theorem~\ref{th:rho_dominance_general} is the first result that provides an analytic, correlation between the control theoretical analysis and real-time implementation.
Primarily, it implies that the constraint dominance from Definition~\ref{def:domination} also carries on to the JSR, giving us a notion of \emph{JSR dominance}.
The results of Theorem~\ref{th:rho_dominance_general} are strategy-independent, 
further reducing the coupling between the control analysis and real-time implementation, and are also independent of the controlled system's dynamics.

Two Corollaries of Theorem~\ref{th:rho_dominance_general} are derived for the commonly used models $\overbar{\genfrac{<}{>}{0pt}{}{m}{\!\:\!\:k\!\:\!\:}}^{\strat}$ and $\overbar{\binom{m}{k}}^{\strat}$, highlighting some practical relations between such constraints.
\begin{corollary}[$\overbar{\binom{m}{k}}^{\strat}$ dominance]
    \label{cor:rho_dominance_mk}
    Given $\lambda^\strat_1 = \overbar{\binom{m}{k_1}}^{\strat}$ and $\lambda^\strat_2 = \overbar{\binom{m}{k_2}}^{\strat}$, if $k_1 \leq k_2$ then
    $
        \rho\bigl(\Alifted_{\lambda^\strat_2}\bigr) \leq \rho\bigl(\Alifted_{\lambda^\strat_1}\bigr).
    $
\end{corollary}
\begin{corollary}[$\overbar{\genfrac{<}{>}{0pt}{}{m}{\!\:\!\:k\!\:\!\:}}^{\strat}$ dominance]
    \label{cor:rho_dominance_cons}
    Given $\lambda^\strat_1 = \overbar{\genfrac{<}{>}{0pt}{}{m}{\!\:\!\:k\!\:\!\:}}^{\strat}$ and $\lambda^\strat_2 =\overbar{\binom{m}{k}}^{\strat}$, then
    $
        \rho\bigl(\Alifted_{\lambda^\strat_2}\bigr) \leq \rho\bigl(\Alifted_{\lambda^\strat_1}\bigr).
    $
\end{corollary}
The conclusions drawn from Theorem~\ref{th:rho_dominance_general} are theoretical, but its practical applicability lies in the algorithm used to find $\rho^{LB}$ and $\rho^{UB}$, i.e., lower and upper bounds for the JSR value.
Using these bounds we can determine the stability of the corresponding switching systems, as follows:
%
$$
\rho^{LB} \bigl( \Alifted_{\lambda^\strat_2} \bigr) \leq \rho \bigl( \Alifted_{\lambda^\strat_2}
\bigr) \leq \rho \bigl( \Alifted_{\lambda^\strat_1} \bigr) \leq \rho^{UB} \bigl(
\Alifted_{\lambda^\strat_1} \bigr).
$$
Regardless of the algorithm used to find the bounds,
if $\lambda^\strat_2 \preceq \lambda^\strat_1$ and $\rho^{UB} ( \Alifted_{\lambda^\strat_1} ) < 1$, the system under $\lambda^\strat_2$ is switching stable.
A similar relation holds for the lower bound.

Theorem~\ref{th:rho_dominance_general} can be further extended by relating the joint spectral radius of a single constraint to sets of constraints.
\begin{theorem}
    \label{th:rho_dominance_set_general}
    Given an arbitrary \ewhc{} $\lambda^\strat$, it holds that
    $
        \rho\left(\Alifted_{\constraints^\strat}\right) \leq \rho \left( \Alifted_{\lambda^\strat} \right) ,\,\, \forall \constraints^\strat \ni \lambda^\strat.
    $
    \begin{proof}
        For an arbitrary \ewhc{} set $\constraints^\strat = \{\lambda^\strat_1, \dots, \lambda^\strat_N\}$, its satisfaction set is $\satis_\ell \,( \constraints^\strat ) = \bigcap_{i \in \{1,\dots, N\}} \satis_\ell \,( \lambda^\strat_i )$.
        Thus, for any $\lambda^\strat_i \in \constraints^\strat$ it holds that 
        $
            \satis_\ell \,( \constraints^\strat ) \subseteq \satis_\ell \,( \lambda^\strat )
        $.
        If a string $b$ is in $\satis_\ell \,( \Lambda^\strat )$ it also belongs to $\satis_\ell \,( \lambda^\strat )$. 
        The set of all possible $A_{b}$ is thus included in the set of all possible $A_{a},\, a \in \satis_\ell \,( \lambda^\strat )$.
        As a consequence it holds that
        \begin{equation*}
            \max_{b \in \satis_\ell \left( \constraints^\strat \right) } \norm{A_{b}}^{1/\ell} \leq
            \max_{a \in \satis_\ell \left( \lambda^\strat \right) } \norm{A_{a}}^{1/\ell}, \quad
            \forall \ell\in\mathbb{N}^{>}.
        \end{equation*}
        The theorem follows immediately when $\ell\rightarrow \infty$.
    \end{proof}
\end{theorem}

As in Theorem~\ref{th:rho_dominance_general}, the more we restrict the execution pattern of the control task with sets of constraints, the lower its JSR will be.
%
Theorem~\ref{th:rho_dominance_set_general} delivers the practical insight that enforcing tighter \ewhc{} to a stable system will \emph{never} destabilise it, as formally stated in the following corollary.
\begin{corollary}
    \label{cor:rho_dominance_set}
    Given an arbitrary \ewhc{} $\lambda^\strat$, if $\rho \left( \Alifted_{\lambda^\strat} \right) < 1$ then
    $       \rho\left(\Alifted_{\constraints^\strat}\right) < 1 ,\,\, \forall \constraints^\strat \ni \lambda^\strat. $
%
\end{corollary}

\section{Evaluation}
\label{sec:evaluation}
We apply the lifted dynamics model presented in Section~\ref{sec:stability} to a representative plant for the process industry, controlled using a PI-controller, sampled with $T = 0.5$~s:
\begin{equation*}
\begin{array}{l}
\begin{cases}
x_{t+1} &= \begin{bmatrix}
                0.606 & 0.304 & 0.076 \\
                0 & 0.606 & 0.304 \\
                0 & 0 & 0.606 \\
            \end{bmatrix} x_t + \begin{bmatrix}
                0.014 \\
                0.091 \\
                0.394 \\
            \end{bmatrix} u_t \\
y_t &= \begin{bmatrix}1 & 0 & 0\end{bmatrix} x_t
\end{cases}\\
\begin{cases}
            z_{t+1} &= z_t + 0.359 y_t \\
            u_{t+1} &= 0.454 z_t + 0.633 y_t.  
\end{cases}
\end{array}
\end{equation*}
\lcss{We analyse the stability of the control systems subject to different $\overbar{\binom{m}{k}}^{\strat}$ constraints.}
We consider all combinations of strategy (Kill or Skip-Next) and actuator mode (zeroing or holding).
For each combination, we generate the lifted set $\Alifted_{\lambda^\strat}$. Its JSR $\rho \left( \Alifted_{\lambda^\strat} \right)$ is then approximated using three different algorithms.
First, a lower and upper bound of $\rho \left( \Alifted_{\lambda^\strat} \right)$ is computed using the {\tt JSR toolbox}~\cite{vankeerberghen2014jsr}.
Then, an upper bound of the JSR is obtained via SOS relaxations, using both the \emph{dense} and \emph{sparse} algorithm from \href{https://github.com/wangjie212/SparseJSR}{\tt SparseJSR}~\cite{sparsejsr}.

Table~\ref{table:stable} displays our results, acquired on an Intel Core i5-8265U@1.60GHz CPU with 8GB RAM.
Lower and upper bounds are denoted ``LB'' and ``UB''. 
All upper bounds obtained with {\tt JSR toolbox} was found greater than the ones obtained with SOS, thus omitted from the Table.
The symbol ``$-$'' means that the SDP solver runs out of memory.
\lcss{%
The SDP solver in {\tt SparseJSR} uses a second-order method.
Thus, a different solver (utilising a first-order method) could reduce memory usage at the cost of potential accuracy loss.}
Bold values represent stable systems under their corresponding \ewhc{}, strategy, and actuator mode.
Starred values represent stable systems inferred from Corollary~\ref{cor:rho_dominance_mk}.
The {\tt JSR toolbox} provides an accurate lower bound and a coarse upper bound.
In contrast, the dense SOS method finds a better upper bound but takes more time.
We compare the time to run both SOS methods, indicating with ``$\times$'' the speedup factor to obtain the sparse bound w.r.t.~the dense.


\begin{table}[h]
\scriptsize
\setlength{\tabcolsep}{3.7pt}
\renewcommand{\arraystretch}{1.06}
\caption{Results obtained with our controlled example.}\label{table:stable}
\vspace{-3mm}
\rowcolors{2}{gray!25}{white}
\begin{center}
\begin{tabular}{|c|clll|clll|}
\hline
\rowcolor{gray!50}
&{JSR\cite{vankeerberghen2014jsr}} &{Dense} & \multicolumn{2}{c|}{Sparse}
&{JSR\cite{vankeerberghen2014jsr}} &{Dense} & \multicolumn{2}{c|}{Sparse}\\
\rowcolor{gray!50}
\multirow{-2}{*}{$\overbar{\binom{m}{k}}^{\strat}$} 
&{LB} &{UB} &{UB} & \multicolumn{1}{c|}{$\times$}
&{LB} &{UB} &{UB} & \multicolumn{1}{c|}{$\times$}\\
$m,$ $k$ & \multicolumn{4}{c|}{\textbf{Kill and zeroing}} & \multicolumn{4}{c|}{\textbf{Kill and holding}} \\
\hline
1, 2
& 0.960 &  1.070 & 1.070 & 0.86
& 0.926 &  1.029 & 1.029 & 0.83\\
1, 3
& 0.920 & \textbf{0.995} & \textbf{0.995} & 0.83
& 0.894 & \textbf{0.971} & \textbf{0.971} & 0.77\\
1, 4
& 0.890 & \textbf{0.945} & \textbf{0.996} & 1.06
& 0.894 & \textbf{0.957} & 1.025$\mathbf{^*}$ & 1.25\\
1, 5
& 0.890 & \textbf{0.922} & \textbf{0.983} & 1.96
& 0.894 & \textbf{0.948} & 1.008$\mathbf{^*}$ & 2.25\\
1, 6
& 0.890 & \textbf{0.920} & \textbf{0.975} & 4.36
& 0.894 & \textbf{0.942} & \textbf{0.995} & 3.68\\
\hline
2, 3
& 0.983 & 1.124 & 1.124 & 0.67
& 0.956 & 1.085 & 1.085 & 0.80\\
2, 4
& 0.960 & 1.079 & 1.079 & 0.74
& 0.927 & 1.039 & 1.039 & 0.86\\
2, 5
& 0.939 & 1.039 & 1.142 & 2.09
& 0.905 & 1.002 & 1.105 & 1.58\\
2, 6
& 0.920 & 1.007 & 1.096 & 12.3
& 0.903 & \textbf{0.974} & 1.080 & 19.2\\
\hline
$m,$ $k$ & \multicolumn{4}{c|}{\textbf{Skip-Next and zeroing}} & \multicolumn{4}{c|}{\textbf{Skip-Next and holding}} \\
\hline
1, 2
& 0.922 &  \textbf{0.924} & \textbf{0.924} & 5.40
& 0.958 &  \textbf{0.958} & \textbf{0.958} & 4.43\\
1, 3
& 0.898 & \textbf{0.974} & \textbf{0.974} & 10.5
& 0.917 & \textbf{0.988} & \textbf{0.988} & 10.4\\
1, 4
& 0.898 & \textbf{0.963} & \textbf{0.963} & 18.2
& 0.890 & \textbf{0.940} & \textbf{0.940} & 15.9\\
1, 5
& 0.898 & \textbf{0.954} & \textbf{0.954} & 17.6
& 0.890 & \textbf{0.929} & \textbf{0.929} & 20.8\\
1, 6
& 0.898 & \textbf{0.946} & \textbf{0.947} & 20.9
& 0.890 & \textbf{0.927} & \textbf{0.927} & 25.8\\
\hline
2, 3
& 0.953 & 1.034 & 1.039 & 4.45
& 0.982 & 1.070 & 1.076 & 5.91\\
2, 4
& 0.922 & 1.033 & 1.040 & 23.9
& 0.958 & 1.079 & 1.086 & 24.2\\
2, 5
& 0.898 & \textbf{0.999} & 1.005 & 77.8
& 0.937 & 1.038 & 1.043 & 58.1\\
2, 6
& 0.907 &{--$\mathbf{^*}$} & 1.007 & \multicolumn{1}{c|}{--}
& 0.917 &{--} & \textbf{0.991} & \multicolumn{1}{c|}{--}\\
\hline

\end{tabular}
\end{center}
\vspace{-5mm}
\end{table}


All the upper bounds computed by {\tt JSR toolbox} are greater than 1, while all lower bounds are \emph{below} 1, thus we cannot draw any conclusion using the {\tt JSR toolbox}.
For all \ewhc{}, $\overbar{\binom{m}{k}}^{\strat}$ where $m=1$ and $2<k\leq 6$ the SOS upper bounds allow us to infer that the system is stable for all combinations of strategy and actuator mode, and also for $k=2$ under the Skip-Next strategy.
From Theorem~\ref{th:rho_dominance_general}, the stability will hold also for all constraints that are harder to satisfy; in particular, Corollary~\ref{cor:rho_dominance_mk} implies stability for all $\overbar{\binom{m}{k}}^{\strat}$ with $m=1$ and $k>6$.
The speedup ratio is growing when $k$ increases, yielding a particularly high benefit of exploiting sparsity for the Skip-Next strategy with actuation zeroing.

\section{Conclusion}
\label{sec:conclusion}
This paper proposes a switching stability analysis framework \lcss{for LTI systems with arbitrary weakly-hard constraints, extending the weakly-hard model and providing an analytic stability} bound.
\lcss{%
The analysis allows us to assess whether computational errors (present in industrial controllers) affect the stability of the controlled systems.
Future work will focus on the performance loss due to the presence of deadline misses following the extended weakly-hard model.
}


\end{document}